\documentclass[article]{amsart}
\usepackage{graphicx}

\newtheorem{theorem}{Theorem}[section]
\newtheorem{lemma}[theorem]{Lemma}

\newtheorem{proposition}[theorem]{Proposition}

\begin{document}

\title{Unavoidable Parallel Minors of 4-Connected Graphs}
\author{Carolyn Chun, Guoli Ding$^*$, Bogdan Oporowski, Dirk Vertigan}
\address{Mathematics Department\\
Louisiana State University\\
Baton Rouge, Louisiana}
\email{$\{$chchchun,ding$\}$@math.lsu.edu}

\date{\today}
\thanks{$^*$ Supported in part by NSF grants DMS-0556091 and ITR-0326387.}

%

\begin{abstract}
A {\it parallel minor} is obtained from a graph by any sequence of edge contractions and parallel edge deletions.  
We prove that, for any positive integer $k$, every internally $4$-connected graph of sufficiently high order contains a parallel minor isomorphic to a variation of $K_{4,k}$ with a complete graph on the vertices of degree $k$, the $k$-partition triple fan with a complete graph on the vertices of degree $k$, the $k$-spoke double wheel, the $k$-spoke double wheel with axle, the $(2k+1)$-rung M\"obius zigzag ladder, the $(2k)$-rung zigzag ladder, or $K_k$.  
We also find the unavoidable parallel minors of $1$-, $2$-, and $3$-connected graphs.
\end{abstract}

\maketitle

\section{Introduction}
\label{introduction}
In this paper, we will explore unavoidable parallel minors in 1-, 2-, 3-, and internally 4-connected graphs of large order, building on the results of~\cite{bogdan} in the last two cases.  
A graph is {\it internally $c$-connected} if it is $(c-1)$-connected and every separating set of order $(c-1)$ divides the graph into exactly one component and one single vertex.  
Since we only consider vertex connectivity in this paper, we may, without loss of generality, restrict our attention to simple graphs, which are graphs containing no loops or parallel edges.  
\\
\\
\indent We will begin by defining some terms and establishing a convenient notation for use throughout this paper.  
All other graph terminology and notation is defined in~\cite{gt}.  
In particular, we denote a subgraph $H$ of graph $G$ as in \cite{gt} by the notation $H\subseteq G$.  
We say that a graph $M$ is a {\it parallel minor} of a graph $G$, written $M\preceq _{\|} G$, if $M$ is obtained by contracting some edges of $G$, then contracting all loops and deleting multiple edges to yield a simple graph.  
A graph $N$ is a {\it minor} of a graph $G$, written $N\preceq G$, if $N$ is a subgraph of a parallel minor of $G$.  
We use the notation $\Phi(G,N)$ to refer to the set $\{ M\preceq _{\|} G : N\subseteq M \text{ and } \vert N \vert = \vert M \vert \}$, where $\vert N\vert$ is the order of $N$.  
Throughout this paper, note that, in order to ensure that $\Phi (G,N)$ is nonempty, $N$ must contain exactly one component in each component of $G$.  
Since we will only use this notation in the context of a connected minor of a connected graph, we will not need to worry about this qualification.  
Observe that $N$ can be obtained from any member of $\Phi(G,N)$ by deleting edges.  
Conversely, a member of $\Phi(G,N)$ is the graph $N$ with extra edges.
\\
\newline \indent The following statement of a Ramsey theorem will be used several times in this paper.  
This theorem tells us what induced subgraphs to expect in large graphs.  

\begin{theorem}
\label{t1}
There is a function $f_{\ref{t1}}$ such that, for any natural number $k$, every graph with order at least equal to the integer $f_{\ref{t1}}(k)$ contains an induced subgraph isomorphic to $K_k$ or $\overline{K_k}$.
\end{theorem}
We will also use the following two theorems concerning $3$- and internally $4$-connected graphs.  
In these theorems, and throughout this paper, $W_k$ is a wheel with $k$ spokes, $D_k$ is a double-wheel with $k$ spokes at each hub, $M_k$ is a M\"obius zigzag ladder with $2k+1$ rungs, and $Z_k$ is a zigzag ladder will $2k$ rungs (see Figure 3).  
These theorems are the main results of~\cite{bogdan}, and are restated for our purposes as follows.  

\begin{theorem}
\label{t2}
There is a function $f_{\ref{t2}}$ such that, for any integer $k$ exceeding two, every 3-connected graph with order at least equal to the integer $f_{\ref{t2}}(k)$ contains a minor isomorphic to $W_k$, or $K_{3,k}$.
\end{theorem}

\begin{theorem}
\label{t3}
There is a function $f_{\ref{t3}}$ such that, for any integers $q$ and $r$ exceeding three, every internally 4-connected graph with order at least equal to the integer $f_{\ref{t3}}(q,r)$ contains a minor isomorphic to $K_{4,q}$, $D_q$, $M_{r}$, or $Z_{r}$.
\end{theorem}

Note that Theorems~\ref{t2} and~\ref{t3} give the sets of unavoidable minors of large 3- and 4-connected graphs, respectively.  
These theorems are an integral step in our determination of the unavoidable parallel minors of 3 and 4-connected graphs.
\\
\newline \indent The four main results of this paper give the sets of unavoidable parallel minors of 1-, 2-, 3-, and 4-connected graphs.  
The families of graphs that we introduce in the figures for use in Theorems~\ref{2c},~\ref{3c}, and~\ref{4c} are discussed in Section~\ref{graphs}, and referred to throughout this paper.  
Our results build on one another, and may be stated as follows.  

\begin{theorem}
\label{1c}
There is a function $f_{\ref{1c}}$ such that, for any positive integer $k$, every connected graph with order at least equal to the integer $f_{\ref{1c}}(k)$ contains a parallel minor isomorphic to $K_{1,k}$, $C_k$, $P_k$, or $K_k$.
\end{theorem}
	\begin{center}
	\includegraphics[width=3.6 in,height=1.53 in]{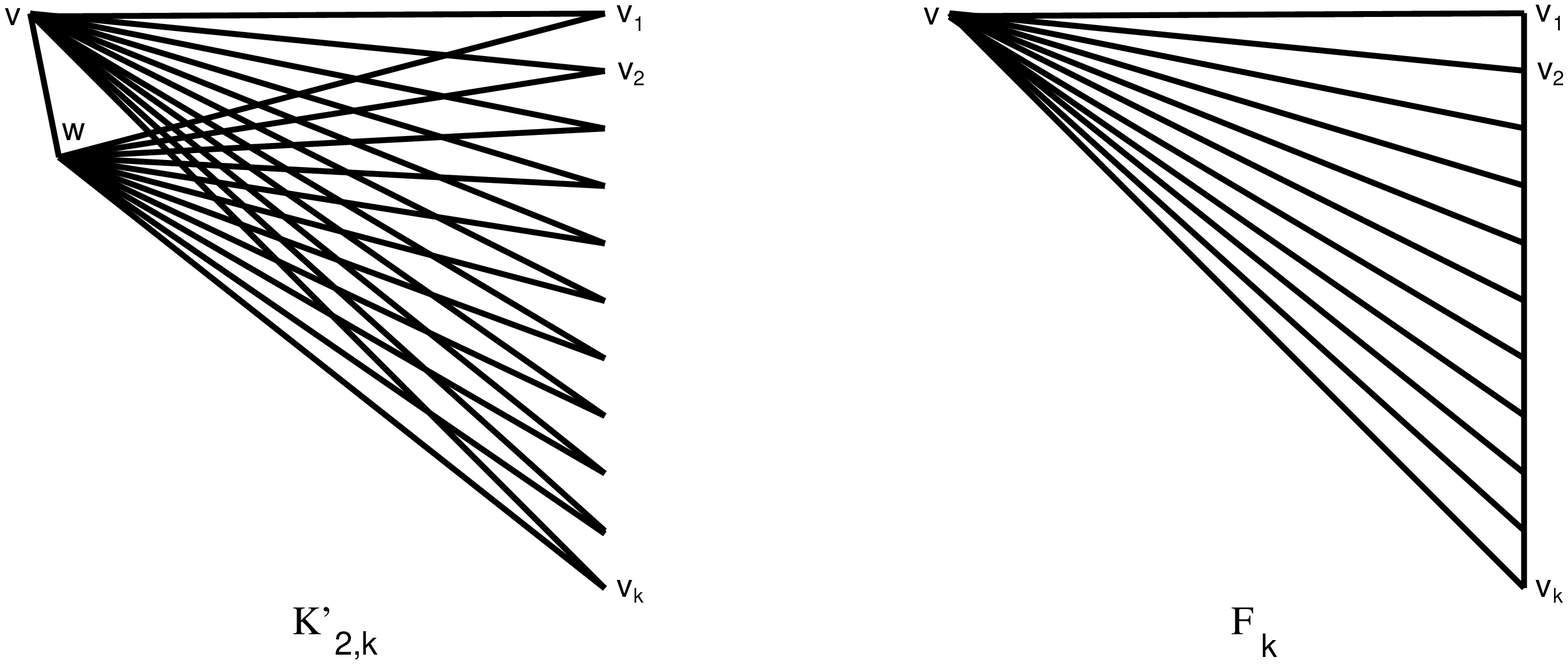}
	
	\begin{scriptsize} FIGURE 1.  Families of 2-connected graphs. \end{scriptsize}
	\end{center}

\begin{theorem}
\label{2c}
There is a function $f_{\ref{2c}}$ such that, for any integer $k$ exceeding two, every $2$-connected graph with order at least equal to the integer $f_{\ref{2c}}(k)$ contains a parallel minor isomorphic to $K_{2,k}'$, $C_k$, $F_k$, or $K_k$.
\end{theorem}
	\begin{center}
	\includegraphics[width=3.6 in,height=1.53 in]{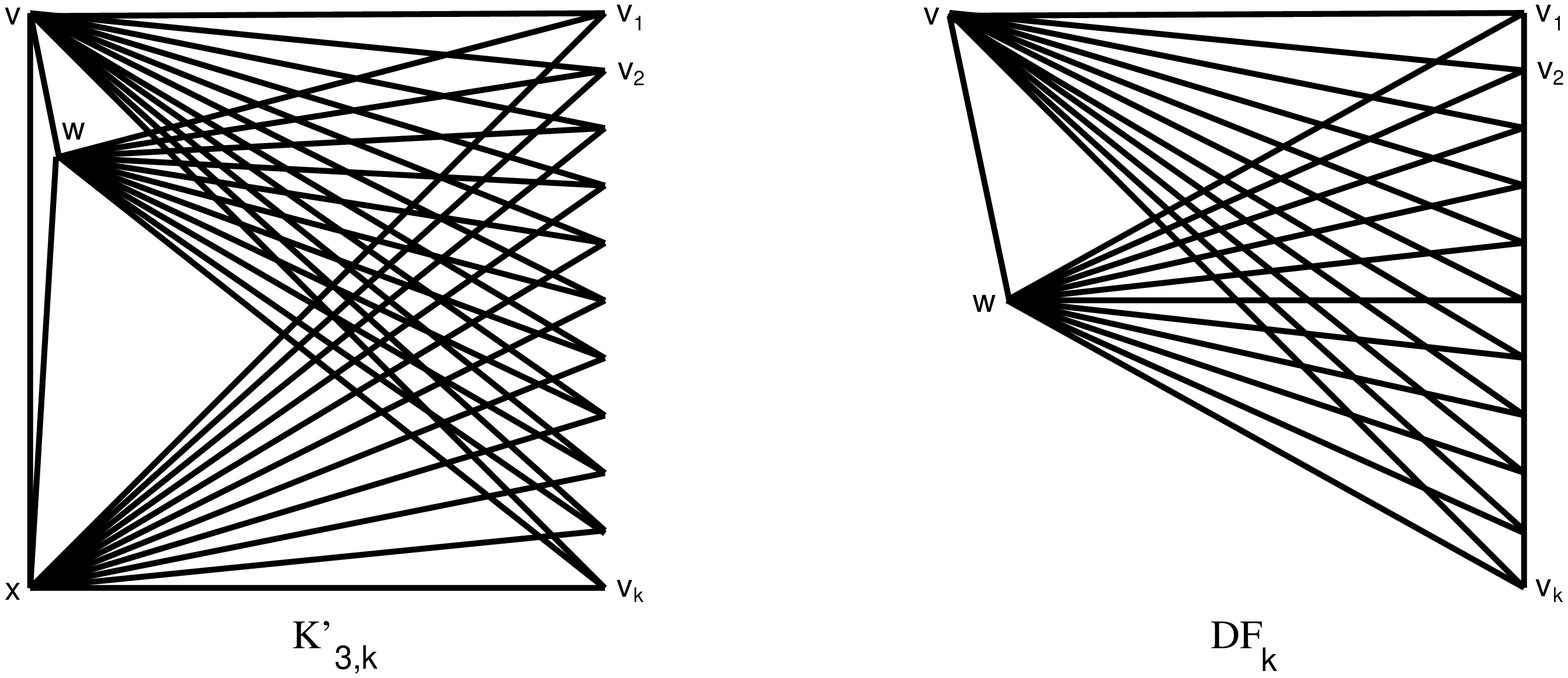}
	
	\begin{scriptsize} FIGURE 2.  Families of 3-connected graphs. \end{scriptsize}
	\end{center}

\begin{theorem}
\label{3c}
There is a function $f_{\ref{3c}}$ such that, for any integer $k$ exceeding three, every $3$-connected graph with order at least equal to the integer $f_{\ref{3c}}(k)$ contains a parallel minor isomorphic to $K_{3,k}'$, $W_k$, $DF_k$, or $K_k$.
\end{theorem}
	\begin{center}
	\includegraphics[width=3.6 in,height=5.4 in]{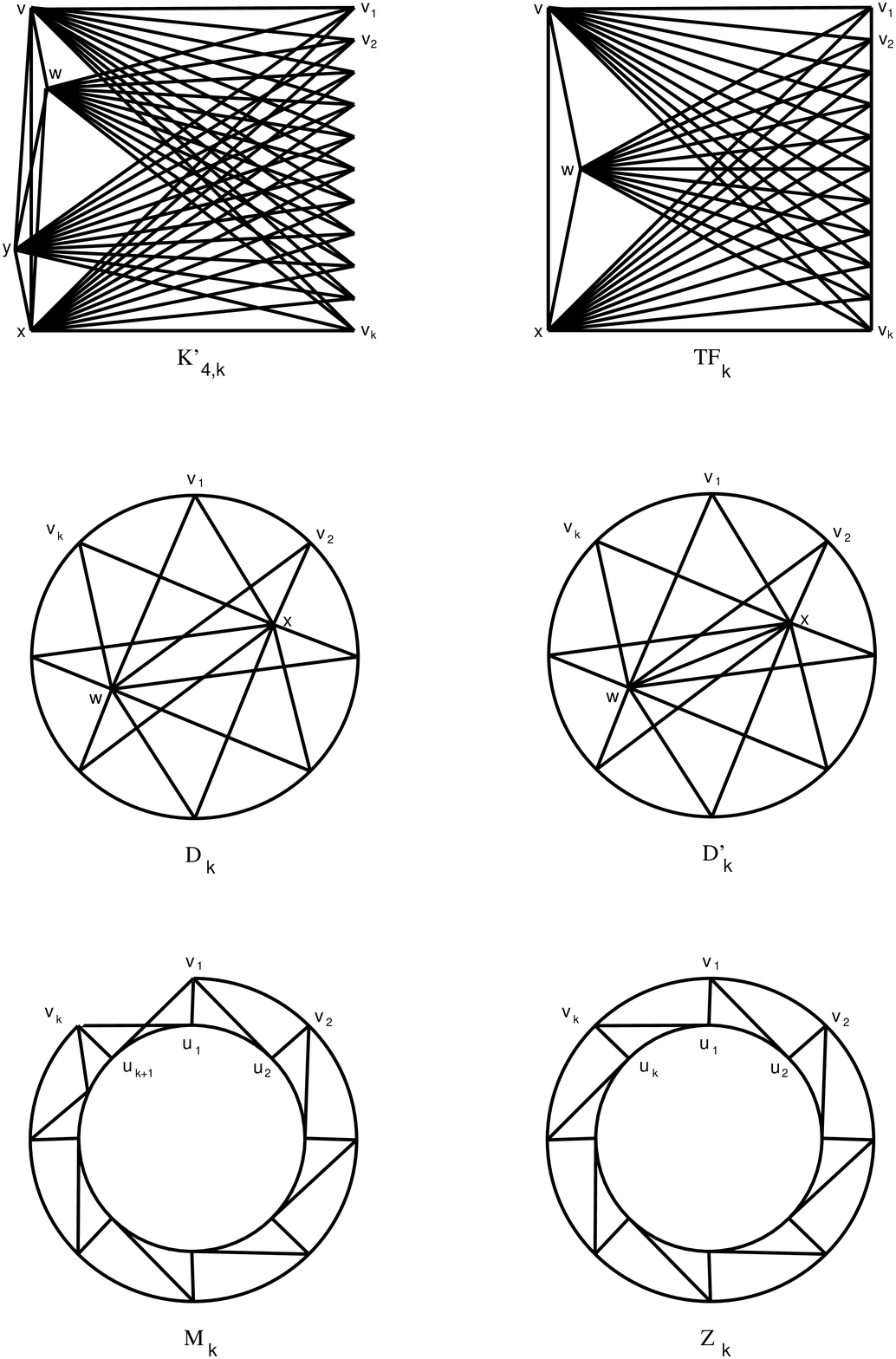}
	
	\begin{scriptsize} FIGURE 3.  Families of 4-connected graphs. \end{scriptsize}
	\end{center}

\begin{theorem}
\label{4c}
There is a function $f_{\ref{4c}}$ such that, for any integer $k$ exceeding four, every internally $4$-connected graph with order at least equal to the integer $f_{\ref{4c}}(k)$ contains a parallel minor isomorphic to $K_{4,k}'$, $D_k$, $D_k'$, $TF_k$, $M_k$, $Z_k$, or $K_k$.
\end{theorem}
\indent Observe that the minors listed in Theorem~\ref{t2} almost form a subset of the minors listed in Theorem ~\ref{3c}, and likewise for Theorem~\ref{t3} and Theorem~\ref{4c}.


\section{Families of 1-, 2-, 3-, and Internally 4-Connected Graphs}
\label{graphs}
We will not prove any theorems in this section:  instead, we will provide a motivation for the specific families of graphs chosen to comprise our sets of unavoidable parallel minors in our variously connected graphs.  
The reader may feel free to turn directly to the internally $4$-connected result, proved in Section~\ref{4proof}, which assumes the $3$-connected result, with the understanding that vertex labeling throughout this paper will follow Figures 1, 2, and 3.
\\
\\
\indent We may have chosen to include the families of graphs from Theorem~\ref{1c} in the list for Theorem~\ref{2c}, since every $2$-connected graph is $1$-connected.  
Observe, however, that each family in the unavoidable set stated in Theorem~\ref{2c} is $2$-connected.  
Likewise, Theorem~\ref{3c} gives a list of families of $3$-connected graphs and Theorem~\ref{4c} gives a list of families of internally $4$-connected graphs.  
We will show each family to be necessary among the unavoidable $c$-connected parallel minors of $c$-connected graphs.
\\
\\
\indent
Consider the $2$-connected family of graphs, $\{ F_k \} _{k>2}$.  
No large parallel minor of a member of this family is $2$-connected, unless it is another member of this same family.  
This is true of each family of graphs listed in Theorem~\ref{2c}.  
Therefore, no family listed contains another in the list.  
The same statement can be made with respect to the $1$-connected graphs listed in Theorem~\ref{1c}, the $3$-connected graphs listed in Theorem~\ref{3c}, and the internally $4$-connected graphs listed in Theorem~\ref{4c}.  
We leave it to the reader to convince himself or herself of this fact.  \\
\\
\indent
We therefore conclude that any set of $1$-connected graphs that comprise an unavoidable set of parallel minors of large, $1$-connected graphs must contain the four families $K_{1,k}$, $C_k$, $P_k$, and $K_k$.  
Similarly, each family in each of the three other sets is necessary.  
The rest of this paper will prove that each set stated in Theorems~\ref{1c}, \ref{2c}, \ref{3c}, and \ref{4c} is sufficient.  
\\
\\
\indent
We make a special remark on Theorem~\ref{4c}.  
Notice that all the listed parallel minors are $4$-connected.  
Since $4$-connectivity implies internal $4$-connectivity, Theorem~\ref{4c} still holds if we replace internal $4$-connectivity with $4$-connectivity.  
In other words, the listed graphs are not only unavoidable in large internally $4$-connected graphs, they are also unavoidable in large $4$-connected graphs.  \\
\\
\indent The unavoidable parallel minors of large, variously connected graphs are significant both because parallel minors are interesting, and because this work complements work done on unavoidable topological minors, the dual operation of parallel minor.  
Unavoidable, topological minors of large, $1$- and $2$-connected graphs are discussed in~\cite{gt}.  
For $3$- and internally $4$-connected graphs, the unavoidable topological minors are described in~\cite{bogdan}.  

\section{Some Parallel Minors in Connected Graphs}
\label{1proof}
In this section, we prove a weak result for $1$-connected graphs, as a step towards proving Theorem~\ref{1c}.  
We begin this section by stating a simple proposition, which provides a set of minors we cannot avoid in a large, connected graph.  
The proof is left to the reader, who may note that the proposition still holds when the word ``minor" is replaced with the word ``subgraph".  
\begin{proposition}
\label{p2}
There is a function $f_{\ref{p2}}$ such that, for positive integers $r$ and $q$, every connected graph of order at least equal to integer $f_{\ref{p2}}(r,q)$ contains $K_{1,r}$, or $P_q$ as a minor.
\end{proposition}
We begin the proof of our connected result by proving the following lemma.  
Theorem~\ref{1c} will be proved in the next section.  
%
%
\begin{lemma}
\label{l1}
There is a function $f_{\ref{l1}}$ such that, for positive integers $k$ and $l$, a connected graph $G$ with order at least equal to integer $f_{\ref{l1}}(k,l)$ contains a parallel minor isomorphic to $K_{1,k}$, $P_k$, or $K_k$; or, $G$ has a $2$-connected graph of order at least $l$ as a parallel minor and has no minor isomorphic to $K_{1,r}$, where $r = f_{~\ref{t1}}(k)$.  
\end{lemma}
\begin{proof}
Let $k$ and $l$ be positive integers.  
We will now select our variables in a particular way to ease the later steps in the proof.  
Let $f_{\ref{t1}}$ and $f_{\ref{p2}}$ be the functions described in Theorem~\ref{t1} and Proposition~\ref{p2} respectively.  
Let $r=f_{\ref{t1}}(k)$, let $q=l(k+1)$, and let $s=f_{\ref{p2}}(r,q)$.  
Set $f_{\ref{l1}}(k,l)=s$.
Let $G$ be a connected graph with order at least $s$.  \\
\\
\indent
Apply Proposition~\ref{p2} to divide our proof into the following two cases, which are exhaustive.  
Case 1:  Graph $G$ contains a minor isomorphic to $K_{1,r}$.  
Case 2:  Graph $G$ contains no minor isomorphic to $K_{1,r}$, and $G$ contains a minor isomorphic to $P_q$.
\\
\\
\indent
If $G$ meets the conditions of Case 1, then take $M\preceq G$ such that $M$ is isomorphic to $K_{1,r}$.  
Fix $H\in \Phi(G,M)$.  
Take vertex $v\in V(H)$ with degree $r$.  
By Theorem~\ref{t1}, the graph $H-v$ has an induced subgraph isomorphic to $K_{k}$ or $\overline{K_{k}}$.\\
\\
\indent
If $H-v$ has an induced subgraph isomorphic to $K_{k}$, then $H$ has a parallel minor isomorphic to $K_k$.  
Assume, therefore, that $H-v$ has an induced subgraph $S$ isomorphic to $\overline{K_{k}}$.  
In $H$, vertex $v$ is adjacent to every other vertex.  
Contract each edge $vu$, where $u\notin V(S)$, deleting the multiple edges and loops, to obtain a parallel minor isomorphic to $K_{1,k}$, which completes Case 1.  
\\
\\
\indent
If $G$ meets the conditions of Case 2, then $G$ has no minor isomorphic to $K_{1,r}$, and we take $M\preceq G$ such that $M$ is isomorphic to $P_q$. 
Fix $H\in \Phi(G,M)$.  
Let $V_\textnormal{c}$ be the set of cut vertices of $H$.  \\
\\
\indent
If $\vert V_\textnormal{c}\vert \geq k+1$, then let $H'$ be obtained recursively from $H$ by contracting, one by one, each edge that is incident with a vertex not in $V_\textnormal{c}$ and deleting loops and multiple edges.  
The parallel minor $H'$ is isomorphic to a path of length at least $k$, hence $G$ has a parallel minor isomorphic to $P_k$.  
We are not finished with Case 2, however, since $H$ may have fewer than $k+1$ cut vertices.  \\
\\
\indent
If $\vert V_\textnormal{c}\vert < k+1$, then there is a large piece of $H$ between cut vertices.  
Let $N$ be a 2-connected subgraph of $H$ of highest order.  
Subgraph $N$ is an end of $H$ or a piece of $H$ between two vertices of $V_\textnormal{c}$, so there are at most $k+1$ places in $H$ to find $N$.  
The order $\vert N\vert$ is therefore at least $\frac{q}{k+1} = l$.  
Let $H'$ be the parallel minor of $G$ obtained from $H$ by contracting, one by one, each edge not in $N$.  
The graph $H'\preceq _{\|} G$ is $2$-connected and has order at least $l$.  
This completes the proof of Lemma~\ref{l1}.
\end{proof}


\section{Unavoidable Parallel Minors of $1$- and $2$-Connected Graphs}
\label{2proof}
We will prove two lemmas before proving the main lemma, Lemma~\ref{l2}, of this section.  
We begin by establishing some notation.  \\
\\
\indent Let $M$ be a minor of $G$, where $M=G/X\backslash Y$.  
Take an edge $e$ in $M$.  
Let $S$ be the set of edges in the multigraph $M\cup Y$ that are parallel with $e$.  
We say that $S$ is the set of edges in $G$ that are in a parallel class with $e$ in $M$.  
\\
\\
\indent For a Hamiltonian parallel minor $M$ of $G$ and a Hamilton cycle $C$ of $M$, the following statements describe an \textit{$H$-set}.  
Let $P$ be a path in $M$ along $C$ such that each vertex of $P$ has degree two in $M$ except for one endpoint of $P$, which may have degree exceeding two.  
Let $e$ be an edge of $P$.  
Let $S$ be a set of edges in $G$ that belong to the same parallel class as edge $e$ in $M$.  
The quintuple $(M, C, S, P,e)$ is an \textit{$H$-set}.  
Furthermore, we say that the $H$-set $(M', C', S', P',e)$ is an \textit{$H$-minor} of the $H$-set $(M, C, S, P,e)$, which we write as $(M',C',S',P',e)\preceq _H(M,C,S,P,e)$, if the following conditions hold.  
\begin{enumerate}
\item $E(C')\subseteq E(C)$.
\item The graph $M'$ is obtained from $M$ by contracting all edges in $E(C - E(C'))$.  
\end{enumerate}
\indent \indent Observe that $C'$ is a Hamilton cycle of $M'$, and that the $H$-minor relation is transitive.  
\\
\\
\indent We say that the \textit{weight} of an $H$-set is the pair $(\vert S\vert , \vert P\vert)$.  
We say that weight $(\vert S\vert , \vert P\vert)$ is greater than weight $(\vert S ' \vert , \vert P ' \vert)$ if $\vert S \vert > \vert S ' \vert$, or if $\vert S \vert = \vert S ' \vert$ and $\vert P \vert > \vert P ' \vert$.  
\\
\\
\indent We will now prove a helpful lemma, which will give us the conditions for finding a longer induced path or a larger parallel set in a Hamiltonian graph by using the $H$-set construction.
\begin{lemma}
\label{sublemma}
For positive integers $d$ and $k$ and for a graph $G$, if $(M, C, S, P,e)$ is an $H$-set such that $\vert M\vert > dk$ and $\Delta (M) < d$, then $\vert P\vert \geq k$, or $H$-set $(M, C, S, P,e)$ has an $H$-minor $(M', C', S', P',e)$ of greater weight such that $\vert M'\vert > \frac{\vert M\vert}{d}$.
\end{lemma}
\begin{proof}
Let $d$ and $k$ be positive integers.  
Let $(M, C, S, P,e)$ be an $H$-set of weight $(\vert S\vert , \vert P\vert)=(\sigma , \pi)$ such that $\vert M\vert=n > dk$, $\Delta (M) < d$, and $\pi < k$.  
By hypothesis, $C$ is the Hamilton cycle of $M$, $e$ is an edge in $P$, which is contained in $C$, and $S$ is a set of $\sigma$ edges in $G$ that are in a parallel class with $e$ in $M$.  
Order the vertices of $C= v_1 v_2 \dots  v_n $ such that the path $P=v_1 v_2 \dots  v_{\pi } $, where $d_M(v_i)=2$ for $i=1,2,\dots , (\pi -1)$.  
Let $e=v_av_{a+1}$.
\\
\\
\indent Consider the neighbors of $v_{\pi}$ in $M$.  
The vertices in $\{ v_{\pi} \} \cup N(v_{\pi})$ divide the cycle $C$ into at most $d$ path segments, since $v_{\pi}$ has fewer than $d$ neighbors.  
There must be a path $v_l v_{l+1} \dots v_{m-1} v_m$ of length at least $\frac{n}{d}>k$ along $C$, with ends in $\{ v_{\pi} \} \cup N(v_{\pi})$ and no other vertices in that set.  
With the following vertex indices, addition is computed modulo $n$.  
\\
\\
\indent In the case where the long path segment contains $P-v_\pi$, index $m$ is equal to $\pi -1$, and we do the following operations.  
Let $M'$ be obtained from $M$ by the contraction of the path $v_{\pi +1} v_{\pi +2} \dots v_{l-1} v_l$ to the vertex $v_{l}$; let $C'$ be the cycle $C$ after this contraction; let $S'= S$; and let $P'=v_1v_2\dots v_{\pi}v_l$.  
The $H$-set $(M', C', S', P',e)$ has weight $(\sigma ,\pi + 1)$ and $(M', C', S', P',e)\preceq _H (M,C,S,P,e)$, which is what we wanted to show.  
\\
\\
\indent We can therefore assume that the long path segment does not meet path $P$.  
In this case, take $f\in E(G)$ such that $f$ is represented by the edge $v_{\pi}v_{m}$.  
Let $S'=S\cup \{ f\}$.  
We obtain $M'$ from $M$ by performing the following contractions.  
\begin{enumerate}
\item Contract the path $v_mv_{(m+1)} \dots v_{a -1} v_{a}$ to vertex $v_a$.  
\item Contract the path $v_{a+1}v_{a+2}\dots v_{\pi}$ to $v_{a+1}$.  
\item Contract the path $v_{\pi +1}v_{\pi +2}\dots v_l$ to vertex $v_l$.   
Note that $l$ is not equal to $\pi$, by construction.  
\end{enumerate}
\indent
\indent Note that vertex $v_{a+1}$ has degree two.  
Let $C'$ be obtained from $C$ by these same contractions, and let $P'= v_av_{a+1}$.  
The $H$-set $(M', C', S', P',e)$ has weight $(\sigma +1,\pi ')$ and $(M', C', S', P',e)\preceq _H (M,C,S,P,e)$, which is what we wanted to show.  
This concludes the proof of Lemma~\ref{sublemma}.
\end{proof}
With the use of this lemma, we will now prove a second lemma.

\begin{lemma}
\label{lemma}
There is a function $f_{\ref{lemma}}$ such that, for integers $k$ and $d$ exceeding two, any Hamilton cycle of a graph with order at least $f_{\ref{lemma}}(k,d)$ contains edges that may be contracted to obtain either a vertex with $d$ neighbors or a parallel minor isomorphic to $C_{k}$.
\end{lemma}
\begin{proof}
Let $k$ and $d$ be integers greater than two.  
Let $r_H= d^{(k-1)(d^2-1)+2}$.  
Set $f_{\ref{lemma}}(k,d)=r_H$.  
Any Hamiltonian graph with at least $r_H$ vertices has a Hamiltonian minor of order $r_H$, so it suffices for our lemma to prove that an arbitrary Hamiltonian graph with order equal to $r_H$ will have our desired structure.  
Let $G_H$ be a Hamiltonian graph of order $r_H$ such that edges of a Hamilton cycle may not be contracted to obtain either a vertex of degree $d$ or a parallel minor isomorphic to $C_{k}$.  
\\
\\
\indent Let $C_H$ be a Hamilton cycle of $G_H$.  
Take vertex $v$ of $C_H$.  
Vertex $v$ has degree less than $d$, so the vertices of $\{ v\} \cup N(v)$ divide $C_H$ into at most $d$ path segments.  
There is some path segment of length at least $\frac{\vert C_H \vert}{d} = \frac{d^{(k-1)(d^2-1)+2}}{d} = d^{(k-1)(d^2-1)+1}$.  
Choose such a path segment, and let $C$ be the cycle obtained from $C_H$ by contracting all edges of $C_H$ that are not in this path segment and that are not incident with $v$.  
Let $G$ be the graph obtained from $G_H$ by the same contractions.  
Observe that $C$ is a Hamilton cycle of $G$, and $\vert G\vert \geq d^{(k-1)(d^2-1)+1}$.  
Without loss of generality, suppose Hamilton graph $G$ to have order exactly $d^{(k-1)(d^2-1)+1}$.  
Let $r=d^{(k-1)(d^2-1)+1}$.  
\\
\\
\indent Let $e$ and $f$ be the two edges in $C$ incident with $v$.  
If $G=C$, then observe that $G$ contains a parallel minor isomorphic to $C_k$.  
We assume not.  
Let $S=\{ e\} $ and let $P$ be the path with endpoint $v$ containing the edge $e$ such that each internal vertex of $P$ has degree two and $P$ has an endpoint with degree exceeding two.  
If $\vert P\vert \geq k$, then we may contract edges in $C-E(P)-\{ e\}$ to obtain a parallel minor isomorphic to $C_{k}$.  
This is forbidden by our assumptions.  
If we find an $H$-set that is an $H$-minor $(M',C',S',P',e)\preceq _H(G,C,S,P,e)$ such that $\vert S' \vert \geq d^2$, then we may contract a path along $C$ in $G$ that contains exactly one end of each edge in $S'$ to obtain a vertex of degree at least $d$.  
This is also forbidden by our assumptions.  
Our restrictions also require that $P'$ have fewer than $k$ vertices, for the same reason that path $P$ does.  
\\
\\
\indent The $H$-set $(G,C,S,P,e)$ has weight at least $(1,1)$, and $\vert P\vert <k$.  
By applying Lemma~\ref{sublemma}, we may find an $H$-set $(M',C',S',P',e)\preceq _H(G,C,S,P,e)$ of greater weight, where $\vert M' \vert > \frac{r}{d} = d^{(k-1)(d^2-1)}$.  
We may do this another $(k-1)(d^2-1)$ times to find a sequence of $H$-sets of strictly increasing weight, each of which is an $H$-minor of the preceding one.  
By our assumptions, for each $H$-set $(M'',C'',S'',P'',e)$ in this sequence, $\vert P''\vert < k$.  
Since this sequence must include at least $(k-1)(d^2-1)+1$ weights greater than $(1,1)$, none of which may repeat, we may apply the pigeonhole principle to conclude that there must be one $H$-set $(M''',C''',S''',P''',e)$ among this sequence with weight greater than $(d^2-1,k-1)$, so that $\vert S''' \vert >d^2-1$.  \\
\\
\indent By transitivity, this $H$-set is an $H$-minor of $(G,C,S,P,e)$.  
This concludes the proof.
\end{proof}

For the final lemma of this section, we will require the following result concerning $2$-connected graphs from ~\cite{gt}, the proof of which is available in the reference.  
The following proposition names two minors which cannot both be absent from a large, 2-connected graph.  
This will provide a natural way of dividing into two cases the 2-connected graphs of high order that we will study in this section.  

\begin{proposition}
\label{p1}
There is a function $f_{\ref{p1}}$ such that, for any integer $r$ exceeding two, every $2$-connected graph of order at least equal to integer $f_{\ref{p1}}(r)$ contains a minor isomorphic to $C_r$ or $K_{2,r}$.
\end{proposition}
As a next step toward proving our 1- and 2-connected results, Theorems~\ref{1c} and~\ref{2c}, we will now prove a lemma concerning 2-connected graphs that is analogous to Lemma~\ref{l1} in the preceding section for connected graphs.

%
%
\begin{lemma}
\label{l2}
There is a function $f_{\ref{l2}}$ such that, for integers $k$ and $q$ exceeding two, every $2$-connected graph with order at least equal to the integer $f_{\ref{l2}}(k,q)$ has a parallel minor isomorphic to $K_{2,k}'$, $C_k$, $F_k$, $K_k$, or a $3$-connected graph of order at least $q$.
\end{lemma}
\begin{proof}
Let $k$ and $q$ be integers exceeding two.  
Let $f_{\ref{l1}}$, $f_{\ref{t1}}$, $f_{\ref{lemma}}$, and $f_{\ref{p1}}$ be the functions described in Lemma~\ref{l1}, Theorem~\ref{t1}, Lemma~\ref{lemma}, and Proposition~\ref{p1} respectively.  
Let $s=f_{\ref{l1}}(k,q)$, $r= f_{\ref{t1}}(k+1) + f_{\ref{lemma}}(k,s)$, and $l=f_{\ref{p1}}(r)$.  
Set $f_{\ref{l2}}(k,q)=l$.  
Let $G$ be a 2-connected graph of order at least $l$.  
\\
\\
\indent
Proposition~\ref{p1} implies that the following two cases are exhaustive.
Case 1:  Graph $G$ has a minor isomorphic to $K_{2,r}$.  
Case 2:  Graph $G$ has no minor isomorphic to $K_{2,r}$, but $G$ has a minor isomorphic to $C_r$.  
\\
\\
\indent
If $G$ meets the conditions of Case 1, then let $M$ be a minor of $G$ that is isomorphic to $K_{2,r}$.  
Fix $H\in \Phi(G,M)$.  
Take $v$ and $w$ in $V(H)$ with degree at least $r$ in $M$.  
By Theorem~\ref{t1}, the graph $H-\{ v,w\}$ has an induced subgraph isomorphic to $K_{k+1}$ or $\overline{K_{k+1}}$.\\
\\
\indent
If $H-\{ v,w\}$ has an induced subgraph isomorphic to $K_{k+1}$, then $H$ has a parallel minor isomorphic to $K_{k}$.  
Assume, therefore, that $H-\{ v,w\}$ has an independent set $X$ of order $(k+1)$.  
In $H$, vertices $v$ and $w$ are adjacent to all other vertices.  
Contract, one by one, each edge that does not have both ends in $X\cup \{ v\}$, deleting the multiple edges and loops, to obtain a parallel minor isomorphic to $K_{2,k+1}$ or $K'_{2,k+1}$.  
Contract any edge to obtain a parallel minor isomorphic to $K'_{2,k}$.  
This completes Case 1.  
\\
\\
\indent
If $G$ meets the conditions of Case 2, then let $M$ be a minor of $G$ that is isomorphic to $C_r$.  
Fix $H\in \Phi(G,M)$.  
The graph $H$ is Hamiltonian.  
\\
\\
\indent Let $C$ be a Hamilton cycle of $H$.  
We may contract edges of $C$ to obtain a parallel minor isomorphic to $C_k$ or a vertex of degree $s$ by Lemma~\ref{lemma}.  
If the former, then we are done, since $C_k$ is among our list of parallel minors.  
If the latter, then contract edges of $C$ to find a vertex of degree $s$.  
This vertex is contained in a Hamiltonian graph, so we can find a minor $N$ of $H$ isomorphic to $F_{s}$.  
Choose $H'\in \Phi(H,N) $.  
\\
\\
\indent
Take vertex $v$ of degree $s$ in $H'$.  
The graph $H'-v$ is connected, so we may apply Lemma~\ref{l1} with the following result.  
The graph $H'-v$ has a parallel minor isomorphic to $K_{1,k}$, $P_k$, $K_k$, or a $2$-connected graph of order at least $q$.  
Therefore, $H'$ has a parallel minor isomorphic to $K'_{2,k}$, $F_k$, $K_k$, or a $3$-connected graph of order at least $q$, respectively.  
This completes Case 2, and the proof of Lemma~\ref{l2}.
\end{proof}
\indent Using Lemma~\ref{l2} with Lemma~\ref{l1}, we will now prove our first major result of this paper, Theorem~\ref{1c}, concerning connected graphs.  
%
%
\begin{proof}[Proof of Theorem~\ref{1c}]
Let $k$ be a positive integer.  
Let $f_{\ref{t1}}$, $f_{\ref{t2}}$, $f_{\ref{l2}}$, and $f_{\ref{l1}}$ be the functions described in Theorem~\ref{t1}, Theorem~\ref{t2}, Lemma~\ref{l2}, and Lemma~\ref{l1} respectively.  
Let $r=f_{\ref{t1}}(k)$, $q=f_{\ref{t2}}(r)$, $l=f_{\ref{l2}}(2k,q)$, and $s=f_{\ref{l1}}(k,l)$.  
Set $f_{\ref{1c}}(k)=s$.  
Let $G$ be a connected graph of order at least $s$.  
\\
\\
\indent
By Lemma~\ref{l1}, graph $G$ has a parallel minor isomorphic to $K_{1,k}$, $P_k$, or $K_k$; or $G$ has a $2$ connected parallel minor of order at least $l$ that has no minor isomorphic to $K_{1,r}$.  
If $G$ has a parallel minor isomorphic to $K_{1,k}$, $P_k$, or $K_k$, then the theorem holds.  
Suppose that $G$ has a $2$-connected parallel minor $H$ of order at least $l$, and $H$ has no minor isomorphic to $K_{1,r}$.  
\\
\\
\indent We apply Lemma~\ref{l2} to $H$ to obtain a $3$-connected parallel minor of $H$ with order $q$, or a parallel minor isomorphic to $K_{2,2k}'$, $C_{2k}$, $F_{2k}$, or $K_{2k}$.  
If $K_{2,2k}'$ is isomorphic to a parallel minor of $H$, then $K_{1,k}$ is isomorphic to a parallel minor of $G$.  
If $C_{2k}$ is isomorphic to a parallel minor of $H$, then $C_k$ is isomorphic to a parallel minor of $G$.  
If $F_{2k}$ is isomorphic to a parallel minor of $H$, then we contract every other spoke of the fan to obtain a parallel minor of $G$ isomorphic to $K_{1,k}$.  
If $K_{2k}$ is isomorphic to a parallel minor of $H$, then $K_k$ is isomorphic to a parallel minor of $G$.  
Therefore, suppose that none of these four parallel minors occur in $G$.  
\\
\\
\indent Let $H'$ be a $3$-connected parallel minor of $H$ with order $q$.  
From Theorem~\ref{t2}, we know that $H$ must have a minor isomorphic to $W_r$ or $K_{3,r}$, so $H$ has a minor isomorphic to $K_{1,r}$, which contradicts our assumption.  
This completes our proof.
\end{proof}
%
%
\indent With the connected result in hand, we continue on to the $2$-connected result.

\begin{proof}[Proof of Theorem~\ref{2c}]
Let $k$ be an integer exceeding two.  
Let $f_{\ref{t2}}$, $f_{\ref{1c}}$ and $f_{\ref{l2}}$ be the functions described in Theorem~\ref{t2}, Theorem~\ref{1c}, and Lemma~\ref{l2}, respectively.  
Let $r=f_{\ref{1c}}(k+2)$, let $q=f_{\ref{t2}}(r)$ and let $l=f_{\ref{l2}}(k,q)$.  
Set $f_{\ref{2c}}(k)=l$.  
Let $G$ be a $2$-connected graph of order at least $l$.  \\
\\
\indent
By Lemma~\ref{l2}, $G$ has a parallel minor isomorphic to $K_{2,k}'$, $C_k$, $F_k$, $K_k$, or a $3$-connected graph of order at least $q$.  
It remains only to investigate the last possibility.  
Let $G$ contain a $3$-connected graph, $G'$, of order at least $q$ as a parallel minor.  
Graph $G'$ has a minor isomorphic to $W_r$ or $K_{3,r}$, by Theorem~\ref{t2}.  \\
\\
\indent
Let $M$ be a minor in $G'$ isomorphic to $W_r$ or $K_{3,r}$.  
Take $H\in \Phi(G,M) $, and take $v\in V(H)$ of degree at least $r$.  
The graph $H$ is $3$-connected, hence $H-v$ is 2-connected.  
Since $H-v$ is connected and $H-v$ has order $f_{\ref{1c}}(k+2)$, the graph $H-v$ has a parallel minor $H'$ isomorphic to $K_{1,k+2}$, $C_{k+2}$, $P_{k+2}$, or $K_{k+2}$, by Theorem~\ref{1c}.  
Since $v$ is non-adjacent to at most two other vertices in $H'$, the graph $H$ must have a parallel minor isomorphic to $K_{2,k}'$, $F_k$, or $K_k$, as required.  
\end{proof}

\section{Unavoidable Parallel Minors of $3$-Connected Graphs}
\label{3proof}
We will now prove our third result, Theorem~\ref{3c}, using our second result, Theorem~\ref{2c}.  
Recall that Theorem~\ref{3c} states that, for an appropriate integer $k$, every $3$-connected graph of high enough order contains $K_{3,k}'$, $W_k$, $DF_k$, or $K_k$ as a parallel minor.  

\begin{proof}[Proof of Theorem~\ref{3c}]
Let $k$ be an integer exceeding three.  
Let $f_{\ref{2c}}$ and $f_{\ref{t2}}$ be the functions described in Theorem~\ref{2c} and Theorem~\ref{t2}, respectively.  
Let $r=f_{\ref{2c}}(k+2)$ and $q=f_{\ref{t2}}(r)$.  
Set $f_{\ref{3c}}(k)=q$.  
Let $G$ be a 3-connected graph of order at least $q$.  
By Theorem~\ref{t2}, the graph $G$ contains a minor $M$ isomorphic to $W_r$ or $K_{3,r}$.  
Choose $H\in \Phi(G,M) $.
\\
\\
\indent
Take $v\in V(H)$ of highest degree.  
Graph $H-v$ is $2$-connected, and has order at least $r$, so $H-v$ contains a parallel minor $H'$ isomorphic to $K_{2,k+2}'$, $C_{k+2}$, $F_{k+2}$, or $K_{k+2}$, by Theorem~\ref{2c}.  
Vertex $v$ is adjacent to all but at most two other vertices in $H'$, hence $G$ has a parallel minor isomorphic to $K_{3,k}'$, $W_k$, $DF_k$, or $K_k$, respectively.  
This completes our proof.
\end{proof}


\section{Unavoidable Parallel Minors of Internally 4-Connected Graphs}
\label{4proof}
Recall that Theorem ~\ref{t3} states the set of unavoidable minors in large, internally 4-connected graphs, $\{ K_{4,q}, D_q, M_{r}, Z_{r} \}$, which will provide the basis for this proof.  
In this section, we will prove our main result, Theorem~\ref{4c}, which states that an internally $4$-connected graph of sufficiently high order will contain as a parallel minor $K_{4,k}'$, $D_k$, $D_k'$, $TF_k$, $M_k$, $Z_k$, or $K_k$.

\begin{proof}[Proof of Theorem ~\ref{4c}]
Let $k$ be an integer exceeding four.  
Let $f_{\ref{3c}}$, $f_{\ref{lemma}}$, and $f_{\ref{t3}}$ be the functions described in Theorem~\ref{3c}, Lemma~\ref{lemma}, and Theorem~\ref{t3}, respectively.  
Let $q=f_{\ref{3c}}(k+3)$, $r=f_{\ref{lemma}}(2k,4q)$, and $n=f_{\ref{t3}}(q,r)$.  
Set $f_{\ref{4c}}(k)=n$.  
Let $G$ be an internally $4$-connected graph of order at least $n$.  
The graph $G$ has a minor isomorphic to $K_{4,q}$, $D_q$, $M_{r}$, or $Z_{r}$, by Theorem~\ref{t3}.
\\
\\
\indent
If $G$ has a minor, $M$, isomorphic to $K_{4,q}$ or $D_q$, then choose $H\in \Phi(G,M)$.  
Take $v$ of highest degree in $H$.  
Graph $H$ is $4$-connected.  
Graph $H-v$ is $3$-connected, so it has a parallel minor $H'$ isomorphic to $K_{3,k+3}'$, $W_{k+3}$, $DF_{k+3}$, or $K_{k+3}$, by Theorem~\ref{3c}.  
Since $v$ is adjacent to all but at most $3$ vertices of $H'$, graph $H$ has a parallel minor isomorphic to $K_{4,k}'$, $D_k$ or $D_k'$, $TF_k$, or $K_k$, respectively.
\\
\\
\indent
Therefore, suppose $G$ has no minor isomorphic to $K_{4,q}$ or $D_q$.  
Then, $G$ has a minor $M$ isomorphic to $M_r$ or $Z_r$.  
These two cases are very similar, so we will present the proof for the case $M \cong Z_r$ and some notes for the $M_r$ case.  
Take $H\in \Phi(G,M)$.  \\
\\
\indent We will work with a collapsed form of $H$.  
Let $H'$ be $H/ \{ v_1u_1, v_2u_2, \dots , v_ru_r \}$ after deleting multiple edges, and let $C$ be the cycle representing the collapsed ladder.  
We apply Lemma~\ref{lemma} to conclude that edges of $C$, a Hamilton cycle, may be contracted to obtain a vertex of degree $4q$ or a parallel minor isomorphic to $C_{2k}$.  \\
\\
\indent
Suppose we can obtain a vertex of degree $4q$ from $H'$ by only contracting edges in $C$.  
Then we may obtain a graph $D$ from $H$ by the contractions of the corresponding pairs of edges in $H$.  
In this case, $D$ contains a vertex $x_1$ of degree at least $2q$.  
Observe that $D$ is a parallel minor of $H$ with the same order as the ladder subgraph contained inside it, so it maintains the ladder structure, which we may label according to Figure~3 with $x$ and $y$ vertices instead of $u$ and $v$ vertices, respectively.  
\\
\\
\indent
The vertex $x_1$ must be adjacent with at least $q$ vertices in either the $x$-vertices or the $y$-vertices of $D$.  
Let $s$ be the order of the $x$-cycle and the $y$-cycle.  
If $x_1$ has $q$ neighbors among the $x$-vertices, then we may contract the path, $y_2y_3 \dots y_{s-2}y_{s-1}$, in $D$ to a vertex $y$ of degree at least $q$.  
The vertices $y$ and $x_1$ are then the two hubs of a minor isomorphic to a subdivision of $D_q$.  
If $x_1$ has $q$ neighbors among the $y$-vertices, then we may contract the path $x_2x_3 \dots x_{s-2}x_{s-1} $ in $D$ to obtain a minor isomorphic to a subdivision of $D_q$.  
We conclude that $D_q\preceq G$, which contradicts our assumptions.  \\
\\
\indent  Suppose we can obtain no vertex of degree $4q$ from $H'$ by contracting edges in the representative Hamilton cycle.  
Then, by Lemma~\ref{lemma}, we must be able to find in $H'$ a parallel minor $N'$ isomorphic to $C_{2k}$ by contracting edges in the Hamilton cycle.  
For every edge that we contract in the Hamilton cycle of $H'$ to obtain $N'$, we contract corresponding pair of edges in $H$ to obtain the graph $N$.  
\\
\\
\indent
Observe that the parallel minor $N\preceq _{ \| } H$ is simply a ladder, possibly with extra edges.  
For convenience, relabel the vertices of $N$ according to the ladder contained inside it (as shown in Figure 3), but with $x$ and $y$ vertices instead of $u$ and $v$ vertices, respectively.  
\\
\\
\indent
Since our original minor $M$ of $G$ is isomorphic to $Z_{r}$, the parallel minor $N$ of $H$ contains a zigzag ladder.  
Since every edge in $N'$ is in the Hamilton cycle of $N'$, the only edges in $N$ which are ``extra,'' ie. not in the ladder, are all edges in $N$ of the form $x_iy_{i+1}$ or $x_{2k}y_1$.  
We can eliminate these interruptions of our zigzag ladder by contracting every other edge in the $x$-cycle and the complementary edges in the $y$-cycle of $N$; that is, contract the edges $u_1u_2$, $u_3u_4$, $\dots$ , $u_{2k-1}u_{2k}$ and the edges $v_2v_3$, $v_4v_5$, $\dots$ , $v_{2k}v_1$.  
In this way we can find a parallel minor of $H$ isomorphic to $Z_k$.  \\
\\
\indent
If $G$ contains a large M\"obius zigzag ladder instead, we contract any triangle in the ladder followed by the remaining rungs to obtain a representative Hamiltonian graph.  
As in the previous case, we will find a cycle in the representative graph corresponding to a large M\"obius ladder that is a parallel minor of $G$, or we will obtain a vertex of high degree by contracting edges of the Hamilton cycle that represents the M\"obius zigzag ladder.  
In the second case, contracting the pairs of edges in the M\"obius ladder corresponding to the edges contracted in the Hamilton cycle will produce a vertex $v$ with high degree in the contracted graph.  
We then contract some long path to obtain a vertex $x$ adjacent with many of the vertices adjacent to $v$, such that these neighbors lie on a cycle that does not contain $v$ or $x$.  
This contracted graph contains a large double wheel, which concludes the argument and completes the proof.
\end{proof}

\section{Related Conclusions and Further Applications}
\label{conclusion}
Since a parallel minor is an induced minor, the reader should note that the set of unavoidable parallel minors in a $c$-connected graph contains the set of unavoidable induced minors.  
With the exception of $C_k$ in the $1$-connected graph case, the families of unavoidable parallel minors and unavoidable $c$-connected induced minors in $c$-connected graphs are identical.  \\
\\
\indent
Parallel minor is the dual matroid operation of series minor.  
The sets of unavoidable series minors in $k$-connected graphs are known for $k\in \{ 1,2,3,4\}$.  
Since any regular matroid can be constructed from graphic and cographic submatroids, together with submatroids isomorphic to $R_{10}$, the results of this paper may contribute to finding the sets of unavoidable minors of variously connected matroids.  \\
\\
\indent
Other avenues of investigation related to the result of this paper include obtaining a set of minors unavoidable in large, $5$-connected graphs.  
Also, it is natural to consider the dual operation of the induced minor, and the unavoidable minors under this dual operation.

\section*{Acknowledgements}
The authors thank James Oxley for helpful discussions concerning this problem.

\end{document}